\documentclass[11pt]{article}
\pdfoutput=1
\usepackage[utf8]{inputenc}
\usepackage{cite}

\usepackage[utf8]{inputenc}
\usepackage[margin=1in]{geometry}
\usepackage[titletoc,title]{appendix}
\usepackage{comment}

\usepackage{amsmath,amsfonts,amssymb,mathtools}
\usepackage{amsthm}

\usepackage{graphicx,float}

\graphicspath{ {./images/} }

\usepackage{csquotes}

\usepackage{amsthm}
\newtheorem{theorem}{Theorem}

\newtheorem{corollary}[theorem]{Corollary}
\newtheorem{lemma}[theorem]{Lemma}
\newtheorem{proposition}[theorem]{Proposition}
\newtheorem{conjecture}[theorem]{Conjecture}
\newtheorem{question}[theorem]{Question}

\theoremstyle{definition}
\newtheorem{definition}{Definition}[section]
\newtheorem{example}{Example}

\theoremstyle{remark}

\usepackage[utf8]{inputenc}
\usepackage[english]{babel}

\newcommand{\R}{\mathbb{R}}
\newcommand{\Z}{\mathbb{Z}}

\newcommand{\tf}{\tilde{F}}
\newcommand{\tfo}{\tilde{f_1}}
\newcommand{\tft}{\tilde{f_2}}

\newcommand{\tg}{\tilde{g}}

\newcommand{\tth}{\tilde{\theta}}

\usepackage{enumerate}

\title{A Family of Universally Optimal Configurations on Rectangular Flat Tori}

\usepackage{authblk}
\author{Nathaniel J. Tenpas}
\affil{Department of Mathematics, Vanderbilt University}
\affil{nathaniel.tenpas@vanderbilt.edu}

\date{}

\begin{document}
\maketitle
\begin{abstract}
    We utilize recently introduced linear programming bounds for the energy of periodic configurations in $\R^d$ to construct configurations which are universally optimal among those of the form $\omega_4+L_\beta$, where $\omega_4$ is a 4-point multiset and $L_\beta\subseteq \R^2$ is a rectangular lattice. For two values of the parameter $\beta$, the configurations are obtained from the $A_2$ lattice. Other values of $\beta$ yield the third, and largest cardinality, class of examples of $\Lambda$-universally optimal configurations which do not arise from lattices $\Lambda$ known or conjectured to be universally optimal in $\R^d$. 
\end{abstract}

\section{Introduction and Overview of Results}
In this paper, we use linear programming bounds to find periodic configurations which are universally optimal among those of the form $\omega_4+L_\beta$, where $L_\beta\subseteq \R^2$ is a rectangular lattice whose aspect ratio depends on the parameter $\beta$. 
 Our use of linear programming bounds and related multivariate polynomial interpolation problems to construct exact minimizers of periodic energy furthers the program begun in \cite{Hardin_2023}.
 
If a lattice 
$\Lambda$ 
is universally optimal among all configurations of fixed density (in the sense of Cohn and Kumar, see \cite{Cohn_Kumar_2007}), then in particular, for any index $n$ sublattice $\Phi$, $\Lambda$ is universally optimal among all configurations of the form $\omega_n+\Phi$. We call this latter condition $\Phi$-universal optimality. For $\beta \in \{\sqrt{3},1/\sqrt{3}\}$, our optimal configuration is either $A_2$ or a scaling and rotation thereof, and thus, our result for these two values of $\beta$ provides new evidence in support of the important conjecture that $A_2$ is universally optimal by showing $A_2$ is $\Phi$-universally optimal for two of its index 4 sublattices.

For other values of $\beta$, our result addresses the question of finding $\Phi$-universally optimal configurations $\omega_n^*+\Phi$ which are not given by some lattice known or conjectured to be universally optimal. Prior to this work, there were only two known classes of examples, both with $n=2$. The configurations obtained from $\beta\not \in \{ \sqrt{3},1/\sqrt{3}\}$ yield a third class of examples with $n=4$, and the first which use linear programming bounds as a significant component of the proof.

The remainder of the introduction is structured as follows: In Section 1.1, we provide the necessary background to state our main theorem. In Section 1.2, we give some context for the result and outline the rest of the paper. 
\subsection{Main Result}
\textit{Periodic energy of finite configurations.} Let   $F:\R^d  \rightarrow (-\infty,\infty]$ be a  lower-semicontinuous potential.  
For a finite multiset $\omega_n=\{x_1,...,x_n\}\subseteq \R^d$ of cardinality $n$, we consider the \textit{$F$-energy}  of $\omega_n$ defined by
\[
E_F(\omega_n):=\sum_{i=1}^n\sum_{\substack{j=1\\ j\neq i}}^n F(x_i-x_j).
\]
 If for some lattice $\Lambda\subseteq \R^d$, $F$ is $\Lambda$-periodic  (i.e., $F(\cdot+v)=F$ for all $v\in \Lambda$), then we also refer to the \textit{$F$-energy} as \textit{periodic energy}.  In this case, without loss of generality, we may assume that $\omega_n$ lies in the flat tori given by some specified fundamental domain $\Omega_\Lambda:=\R^d/\Lambda$, since   replacing a point $x\in\omega_n$ with any point in $x+\Lambda$ does not change $E_F(\omega_n)$.

The \textit{minimal discrete $n$-point $F$-energy} is subsequently defined as
\begin{equation}\label{minEndef}
\mathcal{E}_F(n):= \inf \{E_F(\omega_n)\mid  \omega_n \subseteq \R^d, \, \lvert \omega_n\rvert=n\},
\end{equation}
and an $n$-point configuration $\omega_n$  satisfying $E_F(\omega_n)=\mathcal{E}_F(n)$ is called \textit{$F$-optimal}. Note that lower-semicontinuity and $\Lambda$-periodicity of $F$ along with compactness of $\Omega_\Lambda$ in the flat torus topology imply the existence of at least one $F$-optimal configuration. Periodic energy has been studied in the context of both large $n$ asymptotics 
(cf. \cite{Hardin_2014}, \cite{Borwein_1998})
and exact minimizers for fixed $n$, $\Lambda$, and $F$
 (cf. \cite{Hardin_2023}, \cite{Faulhuber_2023}). We focus on the latter problem in this paper and consider potentials $F$ generated by various $f: [0,\infty)\to [0,\infty]$ with $d$-rapid decay (i.e. $f(r^2) \in \mathcal{O}(r^{-s} ),  r \rightarrow \infty$, for some $s>d$) using 
\begin{equation}\label{Ffdef}
F_{f,\Lambda}(x):=\sum_{v\in \Lambda}f(|x+v|^2).
\end{equation}
Of special importance to us is the case when $f$ is the exponential, mapping $t\rightarrow e^{-a t}$, and we'll denote the resulting potential $F_{a,\Lambda}$.

\textit{Average energy of infinite configurations.}
The periodic energy of a finite configuration $\omega_n$ relates to the average energy of the infinite periodic configuration $\omega_n+ \Lambda$. Following the conventions of \cite{CKMRV_2022},
first take $B(x,r)$ be the ball of radius $r>0$ centered at $x$. If $C$ is an infinite multiset in $\R^d$ such that every ball intersects finitely many points, we call it an \textit{infinite configuration}. Define $C_r:C \cap B(0,r) $ and  the \textit{density of $C$} as
\[
\lim_{r\rightarrow \infty} \frac{\vert C_r \vert}{{\rm Vol}(B(0,r))},
\]
assuming the limit exists and is finite. 
Then for a configuration $C$ of density $\rho$, the \textit{lower $f$-energy of $C$} is 
\[
E^{l}_f(C):=\liminf_{r\rightarrow \infty} \frac{E_f(C_r)}{\vert C_r\vert}.
\] 
If the limit exists, we'll write it as $E_f(C)$ and call it the \textit{average $f$-energy of $C$} or \textit{euclidean energy}. We will also occasionally abuse notation by identifying an $f:\R\rightarrow (-\infty,\infty]$ of $d$-rapid decay with the corresponding radial potential mapping $x\in\R^d$ to $f(\vert x \vert^2)$. If $C'\subseteq \R^d$ has density $\rho$, we say $C'$ is $f$-optimal if $E_f(C')\leq E^{l}_f(C)$ for all configurations $C\subseteq \R^d$ of density $\rho$. In the case that $C$ is periodic of the form $\omega_n+\Lambda$, $C$ approximates a large crystal whose cells are parallelepiped fundamental domains of $\Lambda$, each containing a translation of $\omega_n$. Then $E_f(C)$ is the average energy per particle, and multiplying by the number of particles yields an estimate of the crystal's energy. Thus, the study of euclidean energy for periodic structures and appropriate $f$ can lend insight into the physical phenomena that molecules tend to arrange in crystalline structures at low temperature. 

We have the following connection, which implies that if $\nu_n$ is $F_{f,\Lambda}$-optimal, then $\nu_n+\Lambda$ has minimal average $f$-energy among all configurations of the form $\omega_n+\Lambda$ 
(cf.    \cite[Lemma 9.1]{Cohn_Kumar_2007} or \cite[Chapter 10]{MEbook}).
\begin{proposition}
\label{periodictoaverage}
Let $C=\omega_n+\Lambda$  and $f$ have $d$-rapid decay. 
Then $E_f(C)$ exists and 
\[
E_f(C)= \frac{1}{n}\left(E_{F_{f,\Lambda}} (\omega_{n}^{C})+ n\sum_{0\neq v\in \Lambda} f(v)\right).
\]
\end{proposition}

\textit{Universal optimality.} Much recent study on the energy of infinite configurations has centered around the concept of universal optimality, originally introduced in \cite{Cohn_2003}. 
\begin{definition}
\label{UniOptDef} Let $\Lambda$ be a lattice in $\R^d$, and for each $a>0$, let $f_a:x\rightarrow e^{-a \vert x\vert^2}$. \begin{itemize}
    \item 
We say that an $n$-point configuration $\omega_n\subset \R^d$ is \textit{$\Lambda$-universally optimal} if it is $F_{a,\Lambda}$-optimal for all $a>0$.
\item We say  $\Lambda$ is \textit{universally optimal} if  $\Lambda$ is $f_a$-optimal for all $a>0$.
\end{itemize}
\end{definition}
Both properties are scale and rotation invariant. 
The phrasing of `universal' is apt because a theorem of Bernstein \cite{Bernstein_1929} implies that if $\Lambda$ is universally optimal, then $\Lambda$ minimizes $f$-energy for all $f$ which are completely monotone functions of distance squared with rapid decay and likewise for $\Lambda$-universally optimal configurations. This class of functions includes the Riesz potentials, (mapping $x\rightarrow \vert x\vert^{-s/2}$, $s>d$) for which much is known in both euclidean and compact settings (cf. \cite{Thompson_1897}, \cite{Borodachov_2006}, \cite{Lewin_2022} to name just a few). 
Note by Proposition \ref{periodictoaverage} that a periodic configuration $\omega_n^*+\Lambda$ is  universally optimal among all configurations of the form $\omega_n+\Lambda$ if and only if $\omega_n^*$ is $\Lambda$-universally optimal, and so we will call such a periodic configuration $\Lambda$-universally optimal. How do universal optimality and $\Lambda$-universal optimality relate? If a lattice $\Lambda$ is universally optimal, then in particular it is $\Phi$-universally optimal for each of its sublattices $\Phi$. In fact, using classical physics techniques of \cite{Fisher_1964}, it can be shown that the converse also holds (cf. \cite{Hardin_2023}):
\begin{proposition}
\label{versionsofunivopt}
$\Lambda$ is \textit{universally optimal} if and only if for any of its sublattices 
$\Phi$, we have $\Lambda$ (and equivalently, $\Lambda \cap \Omega_\Phi$) is $\Phi$-universally optimal.
\end{proposition}
Thus, a universally optimal lattice yields infinitely many periodic energy problems for which we can find exact minimizers, and on the other hand, solutions to a particular sequence of periodic energy problems prove the universal optimality of a lattice. \\
Intuitively, the property of universal optimality should be very hard for a lattice to achieve and in certain dimensions ($3,5,6,7$) universally optimal lattices have been shown not to exist. 
In \cite{Cohn_Kumar_2007}, $\Z$ was proved to be universally optimal. More recently, the $E_8$ and Leech lattices were shown to be universally optimal in dimensions 8 and 24 in \cite{CKMRV_2022}. By taking $a\rightarrow \infty$, universal optimality of $\Lambda\subseteq \R^d$ implies that $\Lambda$ solves the $d$-dimensional sphere packing problem, so the results (and methods) in \cite{CKMRV_2022} generalize Viazovska's celebrated solution \cite{Viazovska_2017} of the sphere packing problem in dimension 8. 

\textit{$A_2$ and linear programming.}
The universal optimality results in dimensions $1,8$, and $24$ have been proved with a linear programming approach that stems back to earlier bounds for sphere packing \cite{Cohn_2003} and older Delsarte-Yudin energy bounds for spherical codes (cf. \cite[Chapters 5.5 and 10.4]{MEbook}). It was conjectured  in \cite{Cohn_Kumar_2007} that such an approach could prove the universal optimality of the hexagonal lattice $A_2:=\begin{bmatrix}
    1 & 1/2\\
    0 & \sqrt{3}/2
\end{bmatrix}\Z^2\subseteq \R^2$.  
While this important problem remains open, $A_2$ is rigorously known to be optimal in some settings: it solves the sphere packing problem in $\R^2$ \cite{Toth_1940}, is universally optimal among lattices \cite{Montgomery_88} and optimizes a wide class of Lennard-Jones type potentials \cite{Theil_2006}. 

In \cite{Hardin_2023}, the aforementioned linear programming bounds for euclidean energy were reformulated for periodic energy in a way that reduces the production of sharp bounds to solving multivariate polynomial interpolation problems. Those techniques were then used to show that $A_2$ is both $2A_2$-universally optimal and $L'$-universally optimal, where $L'$ is an index 6 rectangular sublattice. By Proposition \ref{versionsofunivopt}, these results would be implied by the universal optimality of $A_2$.

\textit{Other examples of $\Phi$-universally optimal configurations.}
Proposition \ref{versionsofunivopt} also raises the following question:
\begin{question}
\label{thequestion}
Are there $\Phi$-universally optimal configurations not of the form $\Lambda \cap \Omega_\Phi$ for some $\Lambda$ which is either known or conjectured to be universally optimal?
\end{question}
As far as the author is aware, there are only two known examples (up to scaling and rotation), each of cardinality 2:
\begin{example}
Set $\Lambda=A_2$ and $\omega_2=\{(0,0),(1/2,\sqrt{3}/6)\}$ or $\{(0,0),(1,\sqrt{3}/3)\}$, consisting of the origin and either of the wells of $A_2$, yielding a honeycomb structure (see \cite{Baernstein_1997} for the key result on theta functions or \cite{Su_2015} for a proof in the context of energy). 
\end{example}

\begin{example}
\label{rectexample}
    Let $\Lambda=\{r_1\Z\times \cdots r_d\Z\}$ be a rectangular lattice in $\R^d$ and $\omega_2=\{(0,0),1/2(r_1,\dots,r_d)\}$, consisting of the origin and the centroid of $\Lambda$'s canonical fundamental domain.
\end{example}
Our main result is a proof of a new family of 4-point configurations which are $\Lambda$-universally optimal and address both Question \ref{thequestion} and the universal optimality of $A_2$. Let $L_\beta:= \Z \times 2\beta \Z$ with fundamental domain $\Omega_{\beta}=[0,1]\times[0,2\beta]$. Also set $\Phi_\beta=\begin{bmatrix}
1 & 1/2\\
0 & \beta/2
\end{bmatrix}\Z^2,
$
and 
\[
\omega^*_\beta:=\{(0,0),(0,\beta),(1/2,\beta/2),(1/2,3\beta/2)\}=\Phi_\beta \cap \Omega_\beta .
\]
We use the linear programming bounds of \cite{Hardin_2023} to prove:
\begin{theorem}\label{mainthm}
The configuration $\omega^*_\beta$ is $L_\beta$-universally optimal if and only if $\beta\geq 1/\sqrt{3}$.
\end{theorem}

\begin{figure}
    \centering
    \begin{minipage}{0.45\textwidth}
        \centering
        \includegraphics[height=0.9\textwidth]{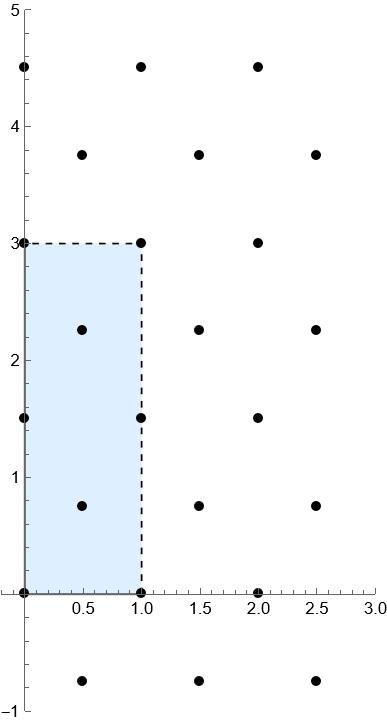} 
        \caption{With $\beta=3/2$ as shown, $\Phi_\beta$ is $L_\beta$-universally optimal.}
    \end{minipage}\hfill
    \begin{minipage}{0.45\textwidth}
        \centering
        \includegraphics[height=0.9\textwidth]{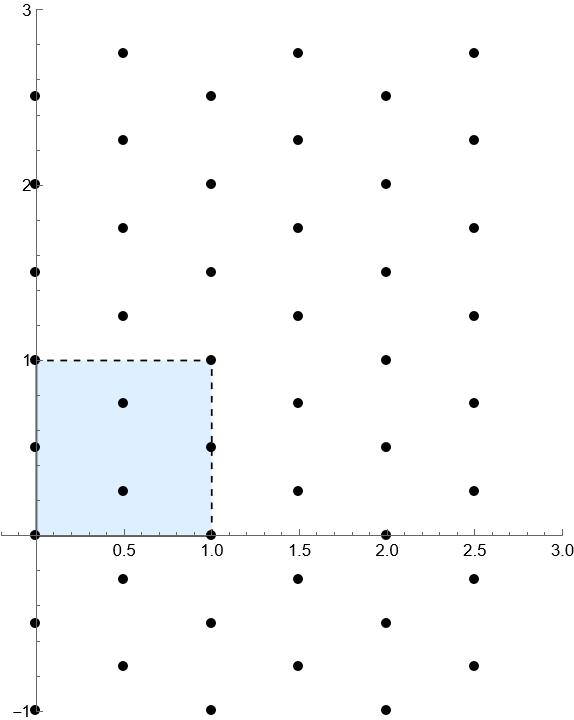} 
        \caption{For $\beta=1/2$, $\Phi_\beta$ is not $L_\beta$ universally optimal.}
    \end{minipage}
\end{figure}

\subsection{Context and Outline}
Note that when $\beta\in \{\sqrt{3},1/\sqrt{3}\}$, $\Phi_\beta$ is either $A_2$
or a scaling and rotation thereof, and so Theorem \ref{mainthm}, like the results of \cite{Hardin_2023}, can be construed as further (on top of already very strong existing) evidence for $A_2$'s universal optimality. In other cases, $\Phi_\beta$ is not a scaling of $A_2$, and so we obtain a third class of examples answering Question \ref{thequestion}. Moreover, our result implies the $\Lambda$-universal optimality of the configurations in Example \ref{rectexample}, as discussed in the final section. 
There, we also present a conjecture that for every $m\geq 2$,
$\Phi_\beta$ yields $\Lambda$-universally optimal $2m$-point configurations for a class of $\Lambda$'s depending on $m$. We'll present some evidence for the conjecture with connections to periodic sphere packing problems.

\textit{Proof outline.}
To prove Theorem \ref{mainthm}, we first introduce the necessary linear programming bounds for periodic energy of \cite{Hardin_2023} in Section \ref{basicssection}, culminating in the construction of a family of bivariate polynomial interpolation problems whose solutions would verify the forward direction of Theorem \ref{mainthm}. In Section \ref{bivinterpprob}, we make an ansatz about the form of potential solutions and then verify that such proposed solutions have the key conditions of conditional positive semi-definiteness and staying below our potential function as required by the linear programming bounds. This verification requires some technical computations which are contained in the appendix. 
In the final section, we discuss results on periodic sphere packing which handle the reverse direction of Theorem \ref{mainthm}. 
\section{Lattices and Linear Programming Bounds}
\label{basicssection}
We begin by introducing background material on lattices in $\R^d$, following the presentation of \cite{Hardin_2023}. Then we present a simplified version of the machinery in \cite{Hardin_2023} where the periodizing lattice is assumed to be rectangular (a more general case is handled in that paper). We will arrive at sufficient conditions, framed as a multivariate polynomial interpolation problem, under which a configuration is $F$-optimal.
\begin{definition}
Let  $\Lambda\subset \R^d$.
  \begin{itemize} 
  \item  $\Lambda$ is a {\em lattice in $\R^d$}  if 
   $\Lambda:= V \mathbb{Z}^d=\left\{\sum_{i=1}^da_iv_i\mid a_1,a_2,\ldots,a_d\in \mathbb{Z}\right\}$ for some nonsingular $d\times d$ matrix $V$  with columns $v_1, \ldots v_d$. We refer to $V$ as a {\em generator} for $\Lambda$. 
   \item Once a choice of generator $V$ is specified, we let  $\Omega_\Lambda:=V\mathbb[0,1)^d$  denote the parallelepiped {\em fundamental domain} for $\Lambda$. The {\em co-volume of $\Lambda$} defined by $\left|\Lambda\right|:=|\det V|$ is the volume  of $\Omega_\Lambda$  and does not depend on our choice of $V$. The lattice $\Phi$ is an index $\kappa$ \textit{sublattice} of $\Lambda$ if $\Phi \subseteq \Lambda$ and $\vert \phi\vert/ \vert \Lambda \vert=\kappa$.
   \item The {\em dual lattice} $\Lambda^*$ of a lattice $\Lambda$ with generator $V$ is the lattice generated by $V^{-T}=(V^T)^{-1}$ or, equivalently,
   $\Lambda^*:=\{v\in \R^d\mid w\cdot v\in \mathbb{Z} \text{ for all } w\in \Lambda    \}.$
   \item We denote by $G^\Lambda$ the group of linear isometries fixing $\Lambda$ and remark that 
   $G^{\Lambda^*}=G^{\Lambda}$.
   \end{itemize}
\end{definition} 

Let $L:=r_1\Z \times \cdots \times r_d \Z$ be a rectangular lattice with each $r_i>0$ and fundamental domain $\Omega_L:=\prod_{i=1}^d [0,r_i)$. Then $H\subseteq G^L$, where $H$ is the order $2^d$ group of coordinate symmetries generated by the maps
\[
h_i(x_1,\dots,x_i,\dots,x_d)=(x_1,\dots, -x_i,\dots, x_d)
\]
for $i=1,\dots,d$.
We have $L^*=\frac{1}{r_1}\Z \times \cdots \times \frac{1}{r_d} \Z$ and define 
\[
W_L=\{ 
(v_1,\dots,v_d)\in L^*: v_i\geq 0, i=1,\dots,d                    \}
\] 
to be a set of representatives of $L^*/H$. 
Continuous functions $g: \R^d\rightarrow (-\infty,\infty]$ which are $L$-periodic are elements (identified with their equivalence classes) of $L^2(\Omega_L)$, the Hilbert space of complex-valued $L$-periodic functions on $\R^d$ with the inner product 
\[
\langle g_1,g_2 \rangle=\frac{1}{r_1\times \cdots \times r_d}\int_{\Omega_L}g_1(x)\overline{g_2(x)}\, dx,
\]
which has the orthogonal basis $\{e^{2\pi i v\cdot x}\mid v\in L^*\}$. Thus, we have the Fourier expansion 
\begin{equation}
 \label{FExp}
g(x)=\sum_{v\in L^*} \hat{g}_v e^{2\pi i v\cdot x}
\end{equation}
with Fourier coefficients $\hat{g}_v:=\langle g,e^{2\pi i v\cdot x} \rangle$
where equality (and the implied unconditional limit on the right hand side) holds in $L^2(\Omega_L)$. If $g$ is further such that $\sum_{v\in L^*}|\hat{g}_v|<\infty$, then the right-hand side of \eqref{FExp} converges uniformly and unconditionally to $g$. Thus, \eqref{FExp} also holds pointwise for every $x\in \R^d$. 

We say that  $g\in L^2(\Omega_\Lambda)$ is {\em conditionally positive semi-definite (CPSD)} if
$\hat{g}_v\ge 0$ for all $v\in\Lambda^*\setminus\{0\}$ and $\sum_{v\in\Lambda^*}\hat{g}_v<\infty$ and say that a CPSD $g$ is {\em positive semi-definite (PSD)} if also $\hat{g}_0\ge 0$.


Given $g\in L^2(\Omega_L)\cap C(\R^d)$, $g$ is $H$-invariant (that is $g(\sigma x)=g(x)$ for all $x$ and $\sigma\in H$) if and only if $\hat g_{\sigma v}=\hat g_{v}$ for all $\sigma \in H$ and $v\in L^*$. The $H$-invariant functions
\begin{align}
C^{H}_v(x):=\frac{1}{|H(v)|}
    \sum_{v'\in H(v)}
    e^{2\pi i   v'\cdot x}, 
    \quad x\in \R^d,
\end{align}
where $H(v)$ is the orbit of $v$ under $H$-translation, allow us to express an arbitrary $g\in L^2(\Omega_L)\cap C(\R^d)$ with $H$-invariance as 

\begin{equation}
    \label{FExpG}
g(x)=\sum_{v\in W_L} |H(v)| \,\hat{g}_v C_v^H(x),
\end{equation}
where equality holds pointwise if $g$ is CPSD. 
\\

We now consider the change of variables
\begin{equation}
\label{littlettransformation}
t_i := \cos(2\pi x_i/r_i), \qquad i=1,...,d, 
\end{equation}
and let $T_{r_1,...,r_d}: \R^d \rightarrow \R^d$ be defined by 
\begin{equation}
\label{Ttransformation}
T_{r_1,\dots,r_d}(x_1,...,x_d):=(t_1,\dots,t_d)=t.
\end{equation}
When context is clear, we refer to this transformation simply as $T$. 
For any $D\subseteq \R^d$, $\tilde{D}:=T(D)\subseteq [-1,1]^d$. 
Similarly, for any $L$-periodic function $h$ with $H$-symmetry, $\tilde{h}$ will refer to the function defined on $[-1,1]^d$ by 
\[
\tilde{h}(t)=h\left(\frac{r_1 \arccos t_1}{2\pi},\dots,\frac{r_d\arccos t_2}{2\pi}\right),
\]
where the assumptions on symmetries of $h$ ensure $\tilde{h}(t)=h(x)$ for all $x\in \R^d$.  We say that $\tilde{h}$ is (C)PSD   if $h$ is (C)PSD.
It is straightforward to see that if $v=k_1/r_1,\dots k_d/r_d\in W_L$, then 
\[
\tilde{C}^{H}_{v}(x)=\prod_{i=1}^d T_{k_i}(t_i),
\]
where $T_j$ is the $j$th Chebyshev Polynomial of the first kind (cf. \cite{Hardin_2023} for a reference).
Furthermore, $\tilde{h}$ is CPSD if and only if its expansion in terms of these polynomials has coefficients that are non-negative (except possibly the $0$th coefficient) and summable. 

Importantly, the map
$F_{a,L}$ behaves nicely under $T$. 
For $c>0$, the classical Jacobi theta function of the third type, is defined by
\begin{equation}
\label{thetasmallaformula}
\theta(c;x):=\sum^{\infty}_{k=-\infty} e^{-\pi k^2 c}e^{2\pi i kx}, \qquad x\in \R.
\end{equation}

Via Poisson Summation on the integers, we have
\begin{equation}
\label{thetabigaformula}
\theta(c;x)=c^{-1/2} \sum^{\infty}_{k=-\infty} e^{-\frac{\pi(k+x)^2}{c}},
\end{equation}

We'll also use 
\[
\tilde{\theta}(c;t):=\theta\left (c,\frac{\arccos t}{2\pi}\right), \quad t\in [-1,1].
\]
It follows from the symmetries of $\theta(c,x)$ that for all $x\in \R$,
\[
\tilde{\theta}(c;\cos 2\pi x)=\theta(c;x),
\]
and moreover, via the Jacobi Triple Product Formula, $\tth$ is known to be absolutely monotone  on $[-1,1]$, as stated below.
\begin{proposition}
\label{absmonotone}
For any $c>0$, the function $\tilde{\theta}=\tilde{\theta}(c; \cdot):[-1,1]\rightarrow(0,\infty)$ is strictly absolutely monotone on $[-1,1]$ and its logarithmic derivative $\tilde{\theta}'/\tilde{\theta}$ is strictly completely monotone on $[-1,1]$.
\end{proposition}

We then have 
\begin{equation}
    \label{Fdecomp}
    \tf_{a,L}(t_1,\dots,t_d)=\prod_{i=1}^d \tth(\pi/(a r_i^2),t_i),
\end{equation}
and so $\tf_{a,L}$ is absolutely monotone in each variable.

Finally, if we have another lattice $\Phi$ such that $L$ is an index $\kappa$ sublattice of $\Phi$, we can consider the $\kappa$-point sets 
$\mu_{\kappa} = \Phi \cap \omega_L$.  
We are now ready to state a simplified set of conditions, applicable when our periodization lattice is rectangular, under which the linear programming bounds in \cite{Hardin_2023} yield optimality of a configuration (see \cite{Hardin_2023}, Cor. 15 for the more general statement).
\begin{corollary}
\label{mukd}
Suppose $L$ is a rectangular lattice with corresponding change of variables  $T:=T_{r_1,\ldots,r_d}$, and $F:\R^d\to (-\infty,\infty]$ is both $H$-invariant and $L$-invariant. Suppose also that
\[
\tg:= c_0+\sum_{v\in W_L} c_v \tilde{C}^{H}_{v}.
\]
is CPSD,
$\tg\leq \tf$ on $[-1,1]^d$, and $L$ is an index $\kappa$ sublattice of some $\Phi$ which satsfies $H\subseteq G^\Phi$.
Then the configuration $\mu_{\kappa}$ 
is $F$-optimal if
\begin{enumerate}
    \item $c_v=0$ for all   $v  \in  \Phi^* \cap W_L$, and 
    \item  $\tg(t) = \tf(t)$ for all $t  \in \tilde \mu_{\kappa}\setminus \{\mathbf{1}\}$ where 
     $\mathbf{1}=(1,1,\ldots,1)\in\R^d$.
\end{enumerate}
\end{corollary}
If such a $\tg$ is constructed, it is called a \textit{magic function}.
We consider the following choices in the setting of $\R^2$: For each $\beta\geq 1/\sqrt{3}$, with 
$L_\beta:= \Z\times 2\beta \Z$ and $\Phi_\beta=\begin{bmatrix}
1 & 1/2\\
0 & \beta/2
\end{bmatrix}\Z^2,
$
as in the statement of Theorem \ref{mainthm}, 
$\omega^*_\beta$ plays the role of $\mu_\kappa$ from Corollary \ref{mukd}. Note that
$H\subseteq G^{\Phi_\beta}$ and 
$\tilde{\omega}^*_\beta/\{\mathbf{1}\}=\{(-1,0),(1,-1)\}$. The maps $\tf$ will be 
\[
\tf_{a,L_\beta}=
\tth(\pi/a_1,t_1)\tth(\pi/(4\beta^2 a_1),t_2)
\] 
for each $\beta$ and $a_1>0$. When $a_1$ and $\beta$ have been fixed, We will define $a_2:=4\beta^2a_1\geq 4a_1/3$,
$\tfo(t_1):=\tth(\pi/ a_1,t_1)$, and $\tft(t_2):=\tth(\pi/a_2,t_2)$ so that 
\[
\tf(t_1,t_2):=\tf_{a,L_\beta}(t_1,t_2)=\tfo(t_1)\tft(t_2).
\]
For each such $\tf$, we will produce a magic $\tg$. To do so, we first compute that 
$\Phi^*=\begin{bmatrix}
    1 &0\\
    1/\beta & 2/\beta\\
\end{bmatrix}\Z^2$,
and thus for $v= (k_1,k_2/(2\beta))\in W_{L_\beta}$, $c_v$ must be zero in our construction of $\tg_a$ only if $2|k_2$ and $k_2/2 \equiv k_1 \pmod 2$.
A priori, these restrictions still leave us with infinitely many degrees of freedom to construct $\tg$
\footnote{The two equality conditions at $(-1,0)$ and $(1,-1)$ and $\tf\leq \tg$ imply that $\frac{\partial (f-g)}{\partial t_2}(-1,0)=0$, and so we need at least a 3-dimensional space to find a magic $\tg$. The author verified that no magic $\tg$ exists in $\text{span}\{ 1,t_1,t_2\}$.}, but we make the ansatz that for every $a_1>0$ and $\beta\geq 1/\sqrt{3}$, a magic $\tg$ exists of the form:
\begin{equation}
    \tg=c_{0,0}+c_{1,0}T_{1}(t_1) +c_{0,1}T_{1}(t_2)+c_{0,2}T_{2}(t_2).
\end{equation}
By proving such a magic $\tg$ exists, we will complete the forward direction of Theorem \ref{mainthm}.
\section{The Bivariate Interpolation Problem}
\label{bivinterpprob}

\subsection{Construction of $\tg$}
Fix $a_1>0$ and $\beta\geq 1/\sqrt{3}$.  It remains to construct a map 
\[
\tg=b_{0,0}+b_{1,0}t_1 +b_{0,1}t_2+b_{0,2}/2(2t_2^2-1)
\]
such that $\tg\leq \tf$ on $[-1,1]$ with equality at $(-1,0)$ and $(1,-1)$ and $b_{1,0},b_{0,1},b_{0,2}\geq 0$.
Our candidate for such a $\tg$ is 
\begin{align}
\tg&=
\tfo'(1)\tft(-1)+\tfo(-1)\tft(0)+
\tfo'(1)\tft(-1)t_1+ 
\tfo(-1)\tft'(0)t_2+
\\
&
t_2^2\left[\tfo(-1)(\tft'(0)-\tft(0))-\tft(-1)(2\tfo'(1)-\tfo(1))\right]\\
&=\tfo(-1)H_{\{-1,0,0\}}[\tft](t_2)+(1+t_1)\tfo'(1)\tft(-1)-t_2^2 \tft(-1)\left(\tfo(-1)-\tfo(1)+2\tfo'(1)\right),
\end{align}
where 
$H_\tau[h](t)$ is the unique Hermite interpolant of degree at most $\vert \tau \vert-1$ to the function $h$ on the node set $\tau$. 
It is straightforward to check that such a $\tg$ satisfies the necessary value conditions at $(-1,0)$ and $(1,-1)$, so it remains to verify the nonegativity of its exapnsion and the inequality $\tf\leq \tg$ on $[-1,1]^2$. 
\subsubsection{Nonnegative Expansion}
Positivity of $b_{0,1}$ and $b_{1,0}$ follow immediately from the absolute monotonicity of $\tfo$ and $\tft$, so it remains to prove  
\begin{lemma}
For any $a_1>0$ and $a_2\geq 4a_2/3$, we have $b_{0,2}>0$.
\end{lemma}
In the appendix, we use our bounds on $\tth$ to construct bounds of 
\[
\frac{\tft'(0)-\tft(0)}{\tft(-1)} \text{ and } \frac{2\tfo'(1)-\tfo(1)}{\tfo(-1)},
\]
which are increasing in $a_1$ and $a_2$ and then use those bounds to show positivity of $b_{0,2}$ via the formula
\[
\frac{b_{0,2}}{\tfo(-1)\tft(-1)}=\frac{\tft'(0)-\tft(0)}{\tft(-1)}-\frac{2\tfo'(1)-\tfo(1)}{\tfo(-1)}.
\]
The positivity of $b_{0,2}$ is equivalent to the fact that for all choices of $a_1>0$ and $b\geq 1/\sqrt{3}$, tangent approximation of $\tf(-1,-1)$ from $(1,-1)$ is more accurate than from $(-1,0)$.

\subsection{$\tf\geq \tg$ on $[-1,1]^2$}
 We begin with just one further technical computation:
\begin{lemma}\label{t1partialat-11}
    For all $a_1>0$ and $a_2\geq 4a_1/3$,
\[
\tfo'(-1)\tft(1)> \tfo'(1)\tft(-1).
\]
\end{lemma}
As with the nonnegativity condition, in the appendix, we prove construct monotone bounds on 
\[
\tfo'(-1)/\tfo'(1) \text{ and } \tft(-1)/\tft(1),
\]
to subsequently prove the lemma. We now may show the desired $\tf\geq \tg$ inequality.
\begin{lemma}\label{staybelowlemma}
For all $a_1>0$, $a_2\geq 4a_1/3$, $(t_1,t_2)\in [-1,1]^2$, we have 
$\tf(t_1,t_2) \geq \tg(t_1,t_2)$.
\end{lemma}

\begin{proof}
We first show that the inequality $\tf\geq \tg$ reduces to a check on the boundary by the second derivative test. Let $H_{\tf}$ denote the Hessian matrix of $\tf$ and likewise for $H_{\tf-\tg}$ with $\tf-\tg$. By the strict complete monotonicity of $\tth$ (see Proposition \ref{absmonotone}), we have 
 $\tilde{f}_i''\tilde{f}_i < (\tilde{f}_i')^2$ for $i\in\{1,2\}$, and thus
\begin{align*}
\det(H_{\tf}&(t_1,t_2)) 
=(\tfo''(t_1) \tfo(t_1))(\tft''(t_2) \tft(t_2))-\tfo'(t_1)^2\tft'(t_2)^2
< 0
\end{align*}
for $t_1,t_2\in[-1,1].$ 
Then we use the positivity of $b_{0,2}$ to compute
\[
\det H_{\tf-\tg}(t_1,t_2)=H_{\tf}(t_1,t_2)-2\tfo''(t_1)\tft(t_2) b_{0,2}<H_{\tf}(t_1,t_2)<0.
\]
Thus, it suffices to verify the difference $\tf-\tg\geq 0$ on the boundary.
When $t_1=-1$, 
 we use the representation
\begin{align}
\tg(-1,t_2)&=\tfo(-1)H_{\{-1,0,0\}}[\tft](t_2)-t_2^2 \tft(-1)\left(\tfo(-1)-\tfo(1)+2\tfo'(1)\right)\\
&\leq \tfo(-1)H_{\{-1,0,0\}}[\tft](t_2)
<\tfo(-1)\tft(t_2),
\end{align}
where the first inequality follows from the fact that $(\tfo(-1)-\tfo(1)+2\tfo'(1))>0$ by the strict absolute monotonicity of $\tfo$, and the second inequality follows from the classical Hermite error formula, which yields
\begin{equation}
\tft(t_2)-H_{\{-1,0,0\}}[\tft](t_2)=\frac{\tft^{3}(\xi)}{3!} (t_2+1)t_2^2\geq 0,
\end{equation}
for some $\xi \in [-1,1]$ depending on $t_2$. In particular, we have $(\tf-\tg)(-1,1)>0$. 
When $t_2=1$, we use the strict absolute monotonicity of $\tfo$ and linearity of $\tg$ in $t_1$ to observe that the difference $(\tf-\tg)(t_1,-1)$ is convex in $t_1$. Thus, the inequalities $(\tf-\tg)(-1,1)>0$ and $\frac{\partial (\tf-\tg)}{\partial t_1}(-1,1)>0$ (from Lemma \ref{t1partialat-11}) immediately yield $(\tf-\tg)(t_1,1)>0$ for all $t_1\in[-1,1]$.
Likewise when $t_2=-1$,  the equalities $(\tf-\tg)(1,-1)=\frac{\partial (\tf-\tg)}{\partial t_1}(1,-1)=0$ immediately yield $(\tf-\tg)(t_1,-1)\geq 0$ for all $t_1\in[-1,1]$.
Finally, we claim that the global minimum of $\tf-\tg$ on $[-1,1]^2$ cannot occur at a point of the form $(1,t_2)$ for $t_2>0$, which would complete the proof. Indeed, for every such point, we have 
\[
\frac{\partial (\tf-\tg)}{\partial t_1}(1,t_2)=\tfo'(1)(\tft(t_2)-\tft(-1))>0
\] 
by strict absolute monoticity of $\tft$, and so $\tf-\tg$ decreases locally as $t_1$ decreases. 
\end{proof}

\section{An Open Problem Related to Sphere Packing}
We conclude with some commentary on our main result, the presentation of an open problem, and the proof of the reverse direction of Theorem \ref{mainthm}.
First, we note that for every $\beta>0$, any $2$-point $\Lambda_\beta$-periodic configuration for $\Lambda_\beta=\begin{bmatrix}
    1&0\\
    0&\beta
\end{bmatrix}\Z^2$ is also a 4-point $L_\beta$-periodic configuration, since $L_\beta$ is an index $2$ sublattice of $\Lambda_\beta$. Since $\omega^*_\beta+L_\beta$ is a 2-point $\Lambda_\beta$-periodic configuration, the universal optimality of $\omega^*_\beta$ among 4-point $L_\beta$-periodic configurations implies its universal optimality among 2-point $\Lambda_\beta$ configurations. Thus, we have given another proof of the universal optimality in Example \ref{rectexample} when $\beta\geq 1/\sqrt{3}$. \\

Recall by Proposition \ref{versionsofunivopt} that the conjectured universal optimality of the $A_2$ lattice would imply that for $\beta=1/\sqrt{3},\sqrt{3}$ and all $m\geq 1$, the $2m$-point configuration $\omega^*_{m,\beta}:=\Phi_\beta \cap ([0,1]\times [0,m\beta])$ is $L^{m}_\beta$-universally optimal, where $L^{m}_\beta:=\begin{bmatrix}
    1&0\\
    0& m\beta
\end{bmatrix}\Z^2$, thus motivating the following conjecture:
\begin{conjecture}
\label{mainconj}
    For all $m\geq 2$, $\omega^*_{m,\beta}$ is   $L^{m}_\beta$-universally optimal if and only if $\beta\geq 1/\sqrt{3}$.
\end{conjecture}
As discussed above, the forward direction of Conjecture \ref{mainconj} holds for all $\beta>0$ when $m=1$. For $m=2$, the conjecture is simply Theorem \ref{mainthm}. Analogous to the euclidean case, for a configuration $\omega_n$ to be $\Lambda$-universally optimal, it must be optimal for the sphere packing problem on the flat torus $\R^d/\Lambda$ (described below). 
A nice result in \cite{Heppes_2000} shows that $\omega^*_{2,\beta}$ is not optimal for sphere packing on the flat torus $\R^2/L^{2}_\beta$ when $\beta < 1/\sqrt{3}$  (in fact, they actually provide the optimal configuration). Thus, Theorem \ref{mainthm} holds. 

\subsection{Connections to Periodic Sphere Packing}
We conclude by giving further evidence for the conjecture in the context of periodic sphere packing. Given a lattice $\Lambda$ and $n$-point configuration $\omega_n=x_1,\dots,x_n$, we define the \textit{$\Lambda$-periodic packing radius} of $\omega_n$ as
\begin{equation}
\label{periodic packing radius}
\delta_\Lambda(\omega_n):=\min \left(\{\vert x_i-x_j+v\vert:i\neq j, v\in \Lambda\} 
\right),
\end{equation}
and the periodic packing problem asks us to find 
\begin{equation}
\label{periodic packing problem}
\mathcal{P}(\Lambda,n):= \sup_{\vert \omega_n\vert=n} \delta_\Lambda(\omega_n)
\end{equation}
and configurations attaining $\mathcal{P}(\Lambda,n)$. We note that that this definition of $\delta_\Lambda(\omega_n)$ 
means that the standard packing radius of the periodic configuration $C=\omega_n+\Lambda$, defined as $\delta(C):=\min_{x\neq y \in C}(\vert x-y \vert )$
satisfies 
\[
\delta(C)=\min \left(\delta_\Lambda(\omega_n),\min_{0\neq v\in \Lambda}(\vert v\vert)\right),
\]
and thus $\delta (C)$ and $\delta_\Lambda(\omega_n)$ may differ.
We choose to omit the vectors $0\neq v\in \Lambda$ from our definition of $\delta_\Lambda(\omega_n)$ 
because they do not depend on the configuration $\omega_n$. 
This problem (sometimes including the vectors $0\neq v\in \Lambda$) has been studied extensively (cf. \cite{Heppes_2000}, \cite{Musin_2015}, \cite{Brandt_2019},\cite{Smirnov_2018} and references therein). 
\\

In \cite{Heppes_2000}, a proof is given for $\beta\in[1/\sqrt{3},\sqrt{3}]$ that the configurations $\delta_{L^{m}_\beta}(\omega^*_{m,\beta})= \mathcal{P}(L^{m}_\beta,2m)$. We give a simiple, alternative proof of the same result for all $\beta\geq 1/\sqrt{3}$ by reducing the dimension of the problem. 
First, using $\beta \geq 1/\sqrt{3}$, 
we check that $\delta_{L^{m}_\beta}(\omega^*_{m,\beta})=\frac{\sqrt{1+\beta^2}}{2}$, achieved for example by the difference of $\vec0$ and $(1/2,\beta/2)$ (this fails for $\beta<1/\sqrt{3}$ as $(0,2\beta)$ is shorter). 
Now take an arbitrary configuration $\omega_{2m}:=x^1,\dots,x^{2m}$ with $x^i=(x^i_1,x^i_2)$ 
for each $i$.
With $d_{\mathbf{T_1}}$ as the flat torus metric from $\R/\Z$,
we have for any $x,y$ that $d_{\mathbf{T_1}}(x,y)\leq \frac12$.
Now let $d_{\mathbf{T_2}}$ be the flat torus metric from $\R/\beta\Z$. The minimal packing radius for $2m$ points in this metric is $\frac{1}{2m}$, achieved only by equally spaced points, so there must be two distinct points 
$x^i,x^j$ such that 
$d_{\mathbf{T_2}} (x^i_2,x^j_2) \leq \frac{\beta}{2}$.
Thus, there exists some $v\in L^{m}_\beta$ such that 
\[
\vert x^i-x^j+v\vert = 
\sqrt{d_{\mathbf{T_1}} (x^i_1,x^j_1)^2+ d_{\mathbf{T_2}} (x^i_2,x^j_2)^2}
\leq \sqrt{(1/2^2)^2+(\beta/2)^2}=\delta_{L^{m}_\beta}(\omega^*_{m,\beta})
\]
as desired.


\appendix
\section{Technical Bounds on Theta functions}
We'll use the following bounds 
developed in the appendix of \cite{Hardin_2023}. For $a\geq \pi^2$, $\epsilon:=1/1000$, $\epsilon_2:=1/50$, $\epsilon_3:=1/40$, we have 
\small
\begin{equation}
\label{thetaboundslargea}
\begin{split}
2e^{-a/4}&< \sqrt{\frac{a}{\pi}}\tth(\pi/a,-1)< 2(1+\epsilon)e^{-a/4},\\
e^{-a/16}&< \sqrt{\frac{a}{\pi}}\tth(\pi/a,0) < (1+\epsilon_2)e^{-a/16},\\
1&< \sqrt{\frac{a}{\pi}}\tth(\pi/a,1) < 1+\epsilon,\\
 \frac{a}{2\pi^2}(a-2)e^{-a/4}
 &< \sqrt{\frac{a}{\pi}}\tth'(\frac{\pi}{a};-1)
 <
 \frac{a}{2\pi^2}(a-2+\epsilon)e^{-a/4},\\
\frac{a e^{-a/16}(1-\epsilon_3)}{4\pi}
&< 
\tth'(\frac{\pi}{a};0)
< 
\frac{a e^{-a/16}}{4\pi},\\
\frac{(1-\epsilon_2)a}{2\pi^2}
 &< \sqrt{\frac{a}{\pi}}\tth'(\frac{\pi}{a};1)
 <
 \frac{a}{2\pi^2}.
\end{split}
\end{equation}
\normalsize
Every bound is directly from the text except the bounds on $\tth(\pi/a,0)$ and $\tth'(\pi/a,0)$, which are straightforward applications of Lemmas 41 and 42 in \cite{Hardin_2023}. 
For $a<\pi^2$, and $d:=\pi^2/a$, we similarly have
\small
\begin{equation}
\label{thetaboundssmalla}
\begin{split}
1-2e^{-d}+99/50e^{-4d}&<\tth(\pi/a,-1)< 1-2e^{-d}+101/50e^{-4d},\\
1-101/50e^{-4d}&<\tth(\pi/a,0)< 1-99/50e^{-4d},\\
1+ 2e^{-d}+99/50e^{-4d}&<\tth(\pi/a,1)< 1+2e^{-d}+101/50e^{-4d},  \\
2e^{-d}-65/8e^{-4d}&<\tth'(\pi/a,-1)<2e^{-d}-63/8e^{-4d},\\
2e^{-d}-1/8e^{-4d}&<\tth'(\pi/a,0)<2e^{-d}+1/8e^{-4d},\\
2e^{-d}+63/8e^{-4d}&<\tth'(\pi/a,1)<2e^{-d}+65/8e^{-4d}.\\
\end{split}
\end{equation}
\normalsize
\section{Positivity of $b_{0,2}$}
To show $b_{0,2}>0$, it suffices to show 
\[\frac{\tft'(0)-\tft(0)}{\tft(-1)}-\frac{2\tfo'(1)-\tfo(1)}{\tfo(-1)}>0.
\]
for all $a_1>0$ and $a_1\geq 4a_2/3$. 
\begin{lemma}
    For $a_1\geq \pi^2$,
\begin{equation}
\label{f1bound1biga}
        \frac{2\tfo'(1)-\tfo(1)}{\tfo(-1)}<
\frac{e^{a_1/4} (a_1-\pi^2)}{2\pi^2}.
\end{equation}
    For $a_2\geq 272\pi/65$, then
\begin{equation}\label{f2bound1bigapos}
\frac{\tft'(0)-\tft(0)}{\tft(-1)}>
\frac{15 e^{3 a_2/16} (65 a_2 - 272 \pi)}{8008 \pi}>0.
\end{equation}
If $\pi^2\leq a_2< 272\pi/65$, then 
\begin{equation}\label{f2bound1biganeg}
\frac{\tft'(0)-\tft(0)}{\tft(-1)}>
 \frac{3 e^{3 a_2/16} (65 a_2 - 272 \pi)}{1600 \pi}.
\end{equation}
Moreover, each of the bounds \eqref{f1bound1biga}, \eqref{f2bound1bigapos}, and \eqref{f2bound1biganeg} is increasing for $a_1$ and $a_2$ at least $\pi^2$, respectively. 
\end{lemma}
\begin{proof}
    The bounds follow immediately from combining the bounds in \ref{thetaboundslargea} and monotonicity comes via  differentiation. For example, 
    \[
\frac{d }{d a_2}\left[
\frac{15 e^{3 a_2/16} (65 a_2 - 272 \pi)}{8008 \pi}
\right]=
\frac{ 15 e^{3 a_2/16} (1040 + 195 a_2 - 816 \pi)}{128128 \pi},
    \]
    whose sign only depends on $1040 + 195 a_2 - 816 \pi$, which is positive for $a_2\geq \pi^2$. The other proofs are similarly straightforward and contained in an accompanying Mathematica notebook \cite{Tenpas_Mathematica_2023}.  
\end{proof}

\begin{lemma}
If $a_1\leq \pi^2$, with $d_1:=\pi^2/a_1$, we have 
\begin{equation}\label{f1bound1smalla}
 \frac{2\tfo'(1)-\tfo(1)}{\tfo(-1)}< \frac{-1 + 1427 e^{-4 d_1}/100 + 2 e^{-d_1}}{1 + 101 e^{-4 d_1}/50 - 2 e^{-d_1}}
\end{equation}
which is increasing in $a_1$. Likewise with $a_2\leq \pi^2$, and $d_2:=\pi^2/a_2$, we have
\begin{equation}\label{f2bound1smalla}
\frac{\tft'(0)-\tft(0)}{\tft(-1)}>\frac{-1 + 371 e^{-4 d_2}/200 + 2 e^{-d_2}}{ 1 + 99 e^{-4 d_2}/50 - 2 e^{-d_2}}
\end{equation}
which is increasing in $a_2$. 
\end{lemma}
\begin{proof}
    The bounds follow immediately from combining the bounds in \ref{thetasmallaformula}, and their monotonicity comes from differentiation. For example, to show the bound in \eqref{f1bound1smalla} is increasing in $a_1$, we show it's decreasing in $d_1$. Its derivative with respect to $d_1$ depends only on the sign of 
\begin{align*}
&(-1 + 1427 e^{-4 d_1}/100 + 2 e^{-d_1})'(1 + 101 e^{-4 d_1}/50 - 2 e^{-d_1})-\\
&(-1 + 1427 e^{-4 d_1}/100 + 2 e^{-d_1})'(1 + 101 e^{-4 d_1}/50 - 2 e^{-d_1})'
\\
&=4887 e^{-5 d_1}/50 - 1629 e^{-4 d_1}/25<0,
\end{align*}
where the last inequality follows simply by checking at $d_1=1$ by Lemma 40 of \cite{Hardin_2023}.
\end{proof}

Now we are ready to show the desired inequality with $b_{0,2}$. If $a_1\geq \pi^2$ (and so also $a_2\geq 4\pi^2/3$, which we note yields the inequality string $a_2\geq 4\pi^2/3>272\pi/65$), then using the fact that $a_2\geq 4a_1/3$, we have the inequalities 
\begin{align*}
\frac{\tft'(0)-\tft(0)}{\tft(-1)}-\frac{2\tfo'(1)-\tfo(1)}{\tfo(-1)}&> 
\frac{15 e^{3 a_2/16} (65 a_2 - 272 \pi)}{8008 \pi}-
\frac{e^{a_1/4} (a_1-\pi^2)}{2\pi^2}\\
&\geq
\frac{15 e^{3 (4a_1/3)/16} (65 (4a_1/3) - 272 \pi)}{8008 \pi}-
\frac{e^{a_1/4} (a_1-\pi^2)}{2\pi^2}\\
&=e^{a_1/4} \frac{-19 \pi^2 + 13 a_1 (-77 + 25 \pi)}{2002 \pi^2}\\
&>0,     
\end{align*}
with the final inequality following by simply by checking $-19 \pi^2 + 13 a_1 (-77 + 25 \pi)>0$ for $a_1=\pi^2$. This point evaluation and others like it are verified in \cite{Tenpas_Mathematica_2023}. 

Now we proceed to an intermediate region, where $a_2 \geq \pi^2$, and $a_1\leq \pi^2$. We have for this region that if $a_2\geq 272\pi/65$, then using monotonicity of the bounds, we obtain 
\begin{align*}
\frac{\tft'(0)-\tft(0)}{\tft(-1)}-\frac{2\tfo'(1)-\tfo(1)}{\tfo(-1)}&> 
0-\frac{-1 + 1427 e^{-4 d_1}/100 + 2 e^{-d_1}}{1 + 101 e^{-4 d_1}/50 - 2 e^{-d_1}}\\
&\geq-\frac{-1 + 1427 e^{-4}/100 + 2 e^{-1}}{1 + 101 e^{-4 }/50 - 2 e^{-1}}>0.
\end{align*}
When $a_2\in [\pi^2, 4\pi^2/3]$, things are more delicate. In that case, if $a_2\in[\pi^2, 272\pi/65]$, since $d_1\geq 4\pi^2/(3a_2)$, we use the inequality
\begin{align*}
\frac{\tft'(0)-\tft(0)}{\tft(-1)}-\frac{2\tfo'(1)-\tfo(1)}{\tfo(-1)}&> \frac{15 e^{3 a_2/16} (65 a_2 - 272 \pi)}{8008 \pi}-\frac{-1 + 1427 e^{-4 d_1}/100 + 2 e^{-d_1}}{1 + 101 e^{-4 d_1}/50 - 2 e^{-d_1}}\\
&\geq 
\frac{15 e^{3 a_2/16} (65 a_2 - 272 \pi)}{8008 \pi}-\frac{-1 + 1427 e^{-4 (4\pi^2/(3a_2))}/100 + 2 e^{-4\pi^2/(3a_2)}}{1 + 101 e^{-4 d_1}/50 - 2 e^{-4\pi^2/(3a_2)}}>0.
\end{align*}
We have a difference of increasing functions in $a_2$, allowing us to verify the final inequality using interval arithmetic in \cite{Tenpas_Mathematica_2023}. \\

Finally, for $a_2\leq \pi^2$, we have $d_1\geq 4d_2/3$, so using monotonoicity of the bounds, we obtain 
\begin{align}
\frac{\tft'(0)-\tft(0)}{\tft(-1)}-\frac{2\tfo'(1)-\tfo(1)}{\tfo(-1)}&>
\frac{-1 + 371 e^{-4 d_2}/200 + 2 e^{-d_2}}{ 1 + 99 e^{-4 d_2}/50 - 2 e^{-d_2}}
-\frac{-1 + 1427 e^{-4 d_1}/100 + 2 e^{-d_1}}{1 + 101 e^{-4 d_1}/50 - 2 e^{-d_1}}\\
&>\frac{-1 + 371 e^{-4 d_2}/200 + 2 e^{-d_2}}{ 1 + 99 e^{-4 d_2}/50 - 2 e^{-d_2}}
-\frac{-1 + 1427 e^{-16 d_2/3}/100 + 2 e^{-4d_2/3}}{1 + 101 e^{-16d_2/3}/50 - 2 e^{-4d_2/3}}
\label{finalline}.
\end{align}
To verify the positivity of \eqref{finalline}, we 
rearrange terms and show 
\begin{align}
    &(-1 + 371 e^{-4 d_2}/200 + 2 e^{-d_2})(1 + 101 e^{-16d_2/3}/50 - 2 e^{-4d_2/3})\\
    &-(-1 + 1427 e^{-16 d_2/3}/100 + 2 e^{-4d_2/3})( 1 + 99 e^{-4 d_2}/50 - 2 e^{-d_2})\\
    &=-(9803/400) e^{-28 d_2/3} + 1629/50 e^{-19 d_2/3} - 599/25 e^{-16 d_2/3} + 
 767 e^{-4 c}/200\label{lemma40use1}\\
 &> e^{-16d_2/3}\left[1529/50 e^{-d_2} - 599/25 + 
 767 e^(4 d_2/3)/200 \right]\label{finalquantity}\\
 &>0
\end{align}
where in \eqref{lemma40use1}, we have used Lemma 40 from \cite{Hardin_2023} to check $-(9803/400) e^{-28 d_2/3} + 2 e^{-19 d_2/3}>0$ for $d_2\geq 1$, and we have checked that $ 1529/50 e^{-d_2} - 599/25 + (
 767 e^{4 d_2/3)}/200 $ is positive for $d_2=1$ and has positive derivative for all $d_2\geq 1$ again using Lemma 40, thus completing the check that $b_{0,2}>0$. 
\section{Proof of Lemma \ref{t1partialat-11}}
We need to show for all $a_1>0$ and $a_2\geq 4a_1/3$ that
\[
\tfo'(-1)/\tfo'(1)> \tft(-1)/\tft(1)
\]
We have the following basic bounds:
\begin{lemma}
For $a_1,a_2\geq \pi^2$, respectively we have
\begin{equation}\label{bound2biga}
\begin{split}
\tfo'(-1)/\tfo'(1)>(a_1-2)e^{-a_1/4}\\
\tft(-1)/\tft(1)< 2(1+\epsilon)e^{-a_2/4}.
\end{split}
\end{equation}
For $a_1,a_2\leq\pi^2$, respectively we have 
\begin{equation}\label{bound2smalla}
\begin{split}
\tfo'(-1)/\tfo'(1)&>\frac{-65/8 e^{-4 d_1} + 2 e^{-d_1}}{65 e^{-4 d_1}/8 + 2 e^{-d_1}},\\
\tft(-1)/\tft(1)&>\frac{1 + 101 e^{-4 d_2}/50 - 2 e^{-d_2}}{1 + 99 e^{-4 d_2}/50 + 2 e^{-d_2}},
\end{split}
\end{equation}
and all of these bounds are decreasing in $a_1$ and $a_2$, respectively.
\end{lemma}
For $a_1\geq \pi^2$, using the monotonicity, we compute
\begin{align*}
    \tfo'(-1)/\tfo'(1)- \tft(-1)/\tft(1)&>
    (a_1-2)e^{-a_1/4}
    -2(1+\epsilon)e^{-a_2/4}\\
    &\geq
    (a_1-2)e^{-a_1/4}
    -2(1+\epsilon)e^{-(4a_1/3)/4}    \\
    &= e^{-a_1/4}(a_1-2-2(1+\epsilon)e^{-a_1/12})\\
    &>e^{-a_1/4}(a_1-4)>0.
\end{align*}
If $a_2\geq \pi^2$ and $a_1\leq \pi^2$,
\begin{align*}
     \tfo'(-1)/\tfo'(1)- \tft(-1)/\tft(1)&>
  \frac{-(65/8) e^{-4 d_1} + 2 e^{-d_1}}{65 e^{-4 d_1}/8 + 2 e^{-d}}
    -2(1+\epsilon)e^{-a_2/4}\\
    &>
\frac{-(65/8) e^{-4}  + 2 e^{-1}}{(65 e^{-4 })/8 + 2 e^{-1}}
    -2(1+\epsilon)e^{-\pi^2/4}>0.
\end{align*}
Finally, for $a_2\leq \pi^2$, we have 
\begin{align}
     \tfo'(-1)/\tfo'(1)- \tft(-1)/\tft(1)&>
  \frac{-(65/8) e^{-4 d_1} + 2 e^{-d_1}}
  {65 e^{-4 d_1}/8 + 2 e^{-d_1}}
  -
  \frac{1 + 101 e^{-4 d_2}/50 - e^{-d_2}}{1 + 99 e^{-4 d_2}/50 + 2 e^{-d_2}}
  \\
    &>
 \frac{-(65/8) e^{-4 (4d_2/3)} + 2 e^{-(4d_2/3)}}
 {65 e^{-4 (4d_2/3)}/8 + 2 e^{-(4d_2/3)}}
  -
  \frac{1 + 101 e^{-4 d_2}/50 - 2 e^{-d_2}}{1 + 99 e^{-4 d_2}/50 + 2 e^{-d_2}}\label{lastone2}.
\end{align}
To prove the positivity of \eqref{lastone2}, as before we rearrange terms to show 
\begin{align*}
&(-(65/8) e^{-4 (4d_2/3)} + 2 e^{-(4d_2/3)})(1 + 99 e^{-4 d_2}/50 + 2 e^{-d_2})\\
&-(1 + 101 e^{-4 d_2}/50 - 2 e^{-d_2})(65 e^{-4 (4d_2/3)} + 2 e^{-(4d_2/3)})\\
&=-(65/2) e^{-28 d_2/3} - 1633/100 e^{-16 d_2/3} + 8 e^{-7 d_2/3}
\\
&>0
\end{align*}
where the last inequality follows from applying Lemma 40 and checking positivity at  $d_2=1$.

\bibliographystyle{abbrv}
\bibliography{citations.bib}

\end{document}